\documentclass[twoside]{report}
\usepackage{amsfonts}

\usepackage{svcon2e}
\usepackage{makeidx}
\usepackage{amsmath}

\newcommand{\ed}{\end{document}}
\makeindex

\input{tcilatex}

\begin{document}

\begin{abstract}
I construct an algebraic model for a typical fiber on a $1+1$ dimensional
spacetime. The vector space comprising the fiber is composed of elements $%
x\otimes x$ formed from the direct product of two copies of an element $x$
in the $D_{2}=C_{2}\otimes C_{2}$ finite group algebra over the real
numbers. The fiber contains subspaces whose elements can be associated with
the tangent and momentum vectors of trajectories in the manifold. The fiber
also contains a subspace whose elements are associated with the local flow
of action of each trajectory. The condition of minimum action translates
into a constraint on the original vector $x$ in the direct product structure.%
\newline
\noindent \textbf{Keywords:} tangent space, manifold, fiber, differential
geometry, $1+1$ space-time.
\end{abstract}

\chapter{Fiber with Intrinsic Action on a\protect\break\ $1+1$ Dimensional
Spacetime}

\chapterauthors{Robert W. Johnson} 
{{\renewcommand{\thefootnote}{\fnsymbol{footnote}} 
\footnotetext{
This work is an outgrowth of a paper presented at the 5th International
Conference on Clifford Algebras and their Applications in Mathematical
Physics, Ixtapa, Mexico, June 27 - July 4, 1999.\newline
AMS Subject Classification: 83C10, 15A66, 53B30.}} }


\section{Introduction}

In a recent paper I described a construction for a vector space with metric.
In this construction one forms elements $x\otimes x$ in a direct product
algebra where $x$ is an element in an underlying finite group algebra \cite
{1}. 
\index{group algebra!finite}One uses a particular decomposition of the
direct product algebra to obtain a direct sum of subspaces. One then
observes that vectors in the various subspaces are interrelated. In
particular examples that I considered there existed a $1d$ subspace whose
measure was determined by the components of a second vector in a higher
dimensional subspace. I proposed that the measure in the $1d$ subspace be
identified with the norm of the vector in the higher dimensional subspace.
In this way both a vector space and its norm are viewed as residing in
particular subspaces of a given direct product algebra $x\otimes x.$ In this
construction the signature of the 
\index{signature!intrinsic} metric arises intrinsically from the particular
underlying finite group algebra.

The realization of Clifford algebras 
\index{Clifford algebra!from finite groups} in terms of underlying finite
group algebras is described by Salingaros \cite{2}. The present work differs
from the approaches of, e.g., Hestenes, Lounesto, and Greider \cite{3} by
considering vector spaces obtained through the decomposition of a direct
product algebra (having elements of form $x\otimes x)$ rather than through a
decomposition of a Clifford algebra. In particular cases that I consider,
however, the underlying finite group algebra (containing elements $x)$ does
correspond to a Clifford algebra.

This work is motivated by the analogous construction in quantum mechanics
where one forms observable vector spaces in terms of bilinear functions of
an underlying state vector. The quantum mechanical $2$-state problem
provides one instructive example \cite{1}. For this problem the underlying
finite group is $C_{4}\otimes H,$ the direct product of the cyclic group 
\index{group!cyclic} 
\index{group!quaternion} of order 4 and the quaternion group. The quantum
wave function $\psi $ is an element in a left ideal of the $C_{4}\otimes H$
group algebra over the real numbers. \footnote{%
This algebra is also termed the complex quaternion algebra. Elements in this
algebra are linear combinations of the $C_{4}\otimes H$ group elements with
real coefficients. The $C_{4}\otimes H$ group multiplication is used along
with the distributive law to induce the algebra product rule.} The
polarization vector and its norm, the total probability, reside in
particular subspaces of the direct product algebra whose elements $\psi
\otimes \psi $ are formed from the product of two copies of the quantum
mechanical state $\psi .$

Vector spaces with metric that are constructed using this procedure are also
of interest as models for typical fiber vector spaces residing at each
location of a configurational manifold such as arises in the context of
classical Lagrangian mechanics. In this paper I develop an algebraic model
for a typical fiber at a location $P$ of a configurational manifold $M$
having 1 space and 1 time dimension. The fiber that I construct contains
both the tangent and momentum vector spaces as subspaces. In addition, a
subspace that can be associated with the flow of action at $P$ is included
and arises in an intrinsic way. The construction for a configuration space
with 2 space dimensions follows in a completely analogous way. The details
for these two 2-dimensional cases are transparent. The extension of this
construction to higher dimensional cases is also straightforward; however,
the multiplicity of subspaces in the corresponding direct product algebras
makes the interpretation of the overall structure more challenging, and it
has not been fully addressed by this author.

In section 2, I review the construction of a vector space with metric
signature $(p,q)=(1,1)$ that corresponds to the tangent space at a point of
a $1+1$ dimensional configurational manifold \cite{1}. This construction
uses the very simple $C_{2}$ group algebra. In section 3 I use the $%
D_{2}=C_{2}\otimes C_{2}$ group algebra to obtain the vector space that is
the primary focus of this paper. In section 4 I summarize these results and
indicate work that still remains to be done.

\section{Algebraic Model for Tangent Space at a Point of a $1+1$ Dimension
Spacetime}

As a starting point for this paper let us review my earlier construction 
\cite{1} of an algebraic model for the vector space with signature $%
(p,q)=(1,1).$ This vector space corresponds to the tangent space 
\index{tangent!space} at each location of a $1+1$ dimensional spacetime $M.$
The tangent space at an arbitrary location $(t_{0},q_{0})\in M$ consists of
tangent vectors 
\index{tangent!vector} to curves passing through $(t_{0},q_{0}).$ Here $%
t_{0} $ and $q_{0}$ denote, respectively, time and spatial location. For a
curve $\gamma =\gamma (\lambda ),\,\lambda \in \mathbb{R},$ the tangent
vector is defined by 
\begin{equation*}
\frac{d\gamma }{d\lambda }=\lim_{\lambda \rightarrow 0}\displaystyle{\frac{%
\gamma (\lambda )-\gamma (0)}{\lambda }}=(\frac{dt}{d\lambda },\frac{dq}{%
d\lambda })
\end{equation*}
where $\gamma (0)=(t_{0},q_{0})$ and $\gamma (\lambda )\in M.$

We begin the construction with the $C_{2}$ group 
\index{group!c@$C_{2}$}which contains two elements $C_{2}=\{\mathbf{1,e}\}$
with $\mathbf{e}^{_{2}}=\mathbf{1},$ and then form the vector space $%
V_{C_{2}}$ whose elements consist of formal sums of the elements of $C_{2}$
with real coefficients. \ An arbitrary element $x\in V_{C_{2}}$ can be
written $x=x_{0}\mathbf{1}+x_{1}\mathbf{e}.$ The product rule, that is
induced by the group multiplication, is used to form the group algebra. We
then consider the direct product of two copies of a vector $x$ in the
algebra $(\mathbf{1} \cong \mathbf{1} \otimes \mathbf{1})$: 
\begin{align*}
x\otimes x & =(x_{0}\mathbf{1}+x_{1}\mathbf{e})\otimes (x_{0}\mathbf{1}+x_{1}%
\mathbf{e}) \\
& =x_{0}^{2}(\mathbf{1\otimes 1})+x_{0}x_{1}(\mathbf{1\otimes e})+x_{0}x_{1}
(\mathbf{e\otimes 1})+x_{1}^{2}(\mathbf{e\otimes e}) \\
& =(x_{0}\mathbf{1}+x_{1}\mathbf{Ee})(x_{0}\mathbf{1}+x_{1}\mathbf{e})
\end{align*}
where the second line follows from bilinearity of the direct product and in
the third line we introduce the notation $\mathbf{E=e\otimes e}$ and $%
\mathbf{e=1\otimes e}.$ $\mathbf{Ee}=\mathbf{e\otimes 1}$ follows from the
product rule in the direct product algebra. Acting on $x\otimes x$ with the
projection operators $P_\pm =\frac12 (\mathbf{1\otimes 1\pm E})$ we obtain 
\begin{equation*}
x\otimes x=[P_{+}(%
\frac{dt}{d\lambda }\mathbf{1}+\frac{dq}{d\lambda }\mathbf{e})+ P_{-}\frac{ds%
}{d\lambda }](\mathbf{1\otimes 1}),
\end{equation*}
where 
\begin{align*}
\frac{dt}{d\lambda }& =x_{0}^{2}+x_{1}^{2} \\
\frac{dq}{d\lambda }& =2x_{0}x_{1} \\
\frac{ds}{d\lambda }& =x_{0}^{2}-x_{1}^{2} = (\frac{dt}{d\lambda }^{2}-\frac{%
dq}{d\lambda }^{2})^{\frac12 }
\end{align*}
The measure $\frac{ds}{d\lambda }$ of the $1d$ $P_{-}x\otimes x$ subspace is
determined up to sign by the two components $(\frac{dt}{d\lambda },\frac{dq}{%
d\lambda })$ of the $2d$ $\ P_{+}x\otimes x$ subspace and can be interpreted
as their norm. We note that the measure $x_{0}^{2}+x_{1}^{2}$ associated
with the $\frac{dt}{d\lambda }$ increment is positive definite and so has an
intrinsic directionality. Also, since $(2x_{0}x_{1})^{2}\leq
(x_{0}^{2}+x_{1}^{2})^{2},$ we have $|\frac{dq}{dt}|\leq 1$ so that there is
a maximum speed.

In this way the product $x\otimes x$ provides a model for an element in the
tangent space of a $1+1$ dimensional spacetime. The set of all such elements 
$x\otimes x$ can be identified with the tangent space itself.

Continuing further we find that rotations of the $2d$ vector space $%
P_{+}x\otimes x$ are induced by acting on $x$ with an element $u=u_{0}%
\mathbf{1}+u_{1}\mathbf{e}$ and forming the product $(x\otimes x)(u\otimes
u)=xu\otimes xu.$ For $u$ such that $u_{0}^{2}-u_{1}^{2}=1,$ $P_{-}x\otimes
x $ is unchanged while the $P_{+}x\otimes x$ vector undergoes a proper
orthochronous rotation.

A completely analogous treatment for the $2d$ Euclidean case is obtained by
substituting the $C_{4}$ group 
\index{group!c@$C_{4}$}algebra for the $C_{2}$ group algebra. This approach
also extends to higher dimensional vector spaces, though in these cases we
encounter a multiplicity of subspaces in the direct product algebra which
make the interpretation of the overall structure more involved \cite{1}.

In the following section we extend this treatment to the case of a typical
fiber at a location $P$ of a configurational manifold $M$ that contains both
the tangent and momentum vector spaces.

\section{Algebraic Model for Fiber on $1+1$ Spacetime with an Intrinsic
Action}

Let us briefly review the classical mechanics that motivate this
construction. The action function 
\begin{equation*}
S_{t_{0},q_{0}}(t,q)=\int_{\gamma }Ldt
\end{equation*}
is the integral of the Lagrangian $L=L(\overset{\cdot }{q},q,t)$ along an
extremal path $\gamma (\lambda ),$ $\lambda \in $ real numbers, connecting
an initial point $(t_{0},q_{0})$ with $(t,q)$ \cite{4}. The action function
for a free particle with mass m in a locally Minkowski coordinate system can
be written as 
\begin{equation*}
S_{t_{0},q_{0}}(t,q)=\int_{\gamma }-m%
\frac{ds}{dt}dt
\end{equation*}
where $ds$ is the proper time with measure $ds=(dt^{2}-dq^{2})^{\frac{1}{2}}$
and $\gamma (\lambda )$ \ is a locally straight line \cite{5}. This action
corresponds to the Lagrangian $L=-m\frac{ds}{dt}$ . The rate of change of
the action along the path $\gamma $ for a fixed initial point is 
\begin{equation}
\frac{dS}{d\lambda }=p\frac{dq}{d\lambda }-H\frac{dt}{d\lambda }
\end{equation}
where 
\begin{equation*}
p=\frac{\partial L}{\partial \overset{\mathbf{\cdot }}{q}}=m\overset{\cdot }{%
q}/(1-\overset{\cdot }{q}^{2})^{\frac{1}{2}}
\end{equation*}
and 
\begin{equation*}
H=p\overset{\mathbf{\cdot }}{q}-L=m/(1-\overset{\mathbf{\cdot }}{q}^{2})^{%
\frac{1}{2}}.
\end{equation*}
We now construct a model for the local tangent and momentum vectors of a
trajectory $\gamma (\lambda )$ on $1+1$ dimensional manifold for which such
an action function arises intrinsically. For this construction we use the
abelian group $(\mathbf{e_{12}}:=\mathbf{e_{1}}\mathbf{e_{2}})$ 
\begin{equation*}
D_{2}=C_{2}\otimes C_{2}=\{\mathbf{1,e}_{\mathbf{1}}\mathbf{,e}_{\mathbf{2}}{%
,\mathbf{e}}_{\mathbf{12}}\}
\end{equation*}
where $\mathbf{e}_{1}^{2}=\mathbf{e}_{2}^{2}=\mathbf{e}_{12}^{2}=\mathbf{1}.$
A general element of the $D_{2}$ group algebra can be written 
\index{group!d@$D_{2}$} 
\begin{equation*}
x=x_{0}\mathbf{1}+x_{1}\mathbf{e}_{\mathbf{1}}+x_{2}\mathbf{e}_{2}+x_{3}%
\mathbf{e}_{12}.
\end{equation*}
\ We decompose $x$ into a sum of two left ideals obtained by acting on the
right with the projection operators $P_{\pm 2}=%
\frac{1}{2}(\mathbf{1}\pm \mathbf{e}_{2}).$ We have 
\begin{equation*}
x=x(P_{+2}+P_{-2})
\end{equation*}
where 
\begin{align*}
xP_{+2}& =[(x_{0}+x_{2})\mathbf{1}+(x_{1}+x_{3})\mathbf{e}_{\mathbf{1}%
}]P_{+2} \\
xP_{-2}& =[(x_{0}-x_{2})\mathbf{1}+(x_{1}-x_{3})\mathbf{e}_{\mathbf{1}%
}]P_{-2}.
\end{align*}
We now form the tensor product of two copies of $x,$ 
\begin{equation*}
x\otimes x=x(P_{+2}+P_{-2})\otimes x(P_{+2}+P_{-2})
\end{equation*}
Expanding out this expression and acting on the left side with the
projection operator $P_{\pm 1}=\frac{1}{2}(\mathbf{1\otimes 1}\pm \mathbf{E}%
_{1}),$ where $\mathbf{E}_{1}$ = $\mathbf{e}_{\mathbf{1}}\otimes \mathbf{e}_{%
\mathbf{1}},$ we obtain, after a change of variables, 
\begin{align*}
x\otimes x=& [P_{+1}(\frac{dt}{d\lambda }\mathbf{1}+\frac{dq}{d\lambda }%
\mathbf{e}_{1})+P_{-1}\frac{ds}{d\lambda }](P_{+2}\otimes P_{+2}) \\[0.5ex]
& +[P_{+1}([\frac{1}{2}(H\frac{dt}{d\lambda }+p\frac{dq}{d\lambda }+m\frac{ds%
}{d\lambda })]^{\frac{1}{2}}\mathbf{1} \\[0.5ex]
& \text{ \ \ \ \ \ \ \ \ }+[\frac{1}{2}(H\frac{dt}{d\lambda }+p\frac{dq}{%
d\lambda }-m\frac{ds}{d\lambda })]^{\frac{1}{2}}\mathbf{e}_{1}) \\[0.5ex]
& +P_{-1}([\frac{1}{2}(H\frac{dt}{d\lambda }-p\frac{dq}{d\lambda }+m\frac{ds%
}{d\lambda })]^{\frac{1}{2}}\mathbf{1} \\[0.5ex]
& \text{ \ \ \ \ \ \ \ \ }-[\frac{1}{2}(H\frac{dt}{d\lambda }-p\frac{dq}{%
d\lambda }-m\frac{ds}{d\lambda })]^{\frac{1}{2}}\mathbf{e}%
_{12})](P_{+2}\otimes P_{-2}+P_{-2}\otimes P_{+2}) \\[0.5ex]
& +[P_{+1}(H\mathbf{1}+p\mathbf{e}_{1})+P_{-1}m](P_{-2}\otimes P_{-2})
\end{align*}
where $\mathbf{e}_{1}=\mathbf{1}\otimes \mathbf{e}_{1}$ and $\mathbf{e}_{12}=%
\mathbf{1}\otimes \mathbf{e}_{12}.$ In this expression 
\begin{align*}
\frac{dt}{d\lambda }& =(x_{0}+x_{2})^{2}+(x_{1}+x_{3})^{2} \\
\frac{dq}{d\lambda }& =2(x_{0}+x_{2})(x_{1}+x_{3}) \\
\frac{ds}{d\lambda }& =(x_{0}+x_{2})^{2}-(x_{1}+x_{3})^{2}=[\frac{dt}{%
d\lambda }^{2}-\frac{dq}{d\lambda }^{2}]^{\frac{1}{2}}
\end{align*}
while for the momenta variables we have, 
\begin{align*}
H& =(x_{0}-x_{2})^{2}+(x_{1}-x_{3})^{2} \\
p& =2(x_{0}-x_{2})(x_{1}-x_{3}) \\
m& =(x_{0}-x_{2})^{2}-(x_{1}-x_{3})^{2}=[H^{2}-p^{2}]^{\frac{1}{2}}
\end{align*}
In analogy to the $C_{2}$ case discussed above, we associate the $%
xP_{+2}\otimes xP_{+2}$ subspace with the tangent to a curve $\gamma
(\lambda )$ on the configurational manifold. The $P_{+1}(xP_{+2}\otimes
xP_{+2})$ portion is identified with the tangent vector $(\frac{dt}{d\lambda 
},\frac{dq}{d\lambda })$ while the $P_{-1}(xP_{+2}\otimes xP_{+2})$ portion
with measure $\frac{ds}{d\lambda }=[\frac{dt}{d\lambda }^{2}-\frac{dq}{%
d\lambda }^{2}]^{\frac{1}{2}}$ is associated with the norm of the tangent
vector.

Similarly, we identify the $xP_{-2}\otimes xP_{-2}$ subspace with the
momentum of the trajectory $\gamma(\lambda)$ at $\lambda.$ The $2d$ $%
P_{+1}(xP_{-2}\otimes xP_{-2})$ projection is identified with the momentum
vector $(H,p)$ while the $1d$ $\ P_{-1}(xP_{-2}\otimes xP_{-2})$ projection
with measure $m=[H^{2}-p^{2}]^{\frac12 }$ is associated with the norm of the
momentum vector.

So far in this development the velocity tangent vector $(\frac{dt}{d\lambda }%
,\frac{dq}{d\lambda})$ is completely independent of the momentum vector $%
(H,p).$ These two vectors contain the 4 degrees of freedom of the original
vector $x$ in the $D_{2}$ algebra.

Let us now consider the $(x\otimes x)(P_{+2}\otimes P_{-2}+P_{-2}\otimes
P_{+2})$ subspace of this algebra. Both the $P_{+1}$and $P_{-1}$ projections
on this subspace are 2-dimensional. Taking the difference of the squares of
the component measures for each of these two 2-vectors we find the resultant 
$m\frac{ds}{d\lambda }.$ The measure squared of the $\mathbf{1}$ component
of the $P_{-1}$ subspace 
\begin{equation*}
\frac12(H\frac{dt}{d\lambda }-p\frac{dq}{d\lambda }+m\frac{ds}{d\lambda })
\end{equation*}
motivates its association with $\frac{dS}{d\lambda }$ of Eq. (1). Minimizing
this measure while keeping $m$ and $\frac{ds}{d\lambda }$ constant and their
product $m\frac{ds}{d\lambda }$ greater than zero requires setting the $%
\mathbf{e}_{12}$ component of the $P_{-1}$ subspace to zero, since the
difference of the squares of the two component measures is fixed. In this
way we obtain the condition 
\begin{equation*}
-m\frac{ds}{d\lambda }=p\frac{dq}{d\lambda }-H\frac{dt}{d\lambda }
\end{equation*}
that corresponds to eq. 1. This equation links the tangent vector $(\frac{dt%
}{d\lambda },\frac{dq}{d\lambda })$ and the momentum vector $(H,p).$ It is
equivalent to the condition $\frac{dq}{dt}=\frac{p}{E}$ that holds for the
trajectory of a free particle. Interestingly, this condition translates to
the requirement 
\begin{equation*}
x=\frac{1}{x_{0}}(x_{0}+x_{1}\mathbf{e}_{1})(x_{0}+x_{2}\mathbf{e}_{2})\;%
\text{ for }\;x_{0}\neq 0
\end{equation*}
for the original vector $x$ in the $D_{2}$ group algebra; or, stated
differently, it ensures that $x$ can be written as the\ tensor product of
two vectors in the $C_{2}$ group algebra.

Transformations that preserve the norms of the tangent and momentum vectors
and the action differential $\frac{dS}{d\lambda}=-m\frac{ds}{d\lambda}$ can
be induced by acting on $x$ with an element $u=u_{0}+u_{1}\mathbf{e}%
_{1}+u_{2} \mathbf{e}_{2}+u_{3}\mathbf{e}_{3}$ in the $D_{2}$ algebra and
forming the product $(x\otimes x)(u\otimes u)=xu\otimes xu.$ For $u$ such
that $u=u_{0}+u_{1} \mathbf{e}_{1}$ and $u_{0}^{2}-u_{1}^{2}=1,$ the $%
P_{-}x\otimes x$ subspace that contains the three norms is unchanged. The
tangent and momentum vectors undergo a proper orthochronous transformation
under this rule.

Finally, we also note that if instead of using $D_{2}=C_{2}\otimes C_{2}$ in
the construction above, we use 
\begin{equation*}
C_{2}\otimes C_{4} = \{\mathbf{1,e}_{2}\}\otimes \{\mathbf{1,e}_{1}\} = \{%
\mathbf{1,e}_{1}\mathbf{,e}_{2}\mathbf{,e}_{12}\}
\end{equation*}
where $\mathbf{e}_{12}=\mathbf{e}_{1}\mathbf{e}_{2},$ $\mathbf{e}_{21}=%
\mathbf{e}_{12},$ $\mathbf{e}_{2}^{2}=+1$ and $\mathbf{e}_{1}^{2}=\mathbf{e}%
_{12}^{2}=-1$ then we obtain the $2d$ Euclidean counterpart to the above $%
1+1 $ spacetime case.

\section{\protect\bigskip Summary}

We associate the vector $x\otimes x$ in the $C_{2}\otimes C_{2}$ group
algebra with an element in a typical fiber residing at a point in a $1+1$
dimensional configurational manifold: 
\begin{align*}
x\otimes x = & \{P_{+1}(\frac{dt}{d\lambda}\mathbf{1}+ \frac{dq}{d\lambda }%
\mathbf{e}_{1})+ P_{-1}\frac{ds}{d\lambda}\} (P_{+2}\otimes P_{+2}) \\[0.5ex]
& +\{P_{+1}{\Big(}(H\frac{dt}{d\lambda })^{\frac12}\mathbf{1} +(p\frac{dq}{%
d\lambda })^{\frac12 }\mathbf{e}_{1}{\Big)} \\[0.5ex]
& +P_{-1}(m\frac{ds}{d\lambda })^{\frac12 }\} (P_{+2}\otimes P_{-2}+
P_{-2}\otimes P_{+2}) \\[0.5ex]
& +\{P_{+1}(H\mathbf{1}+p\mathbf{e}_{1})+P_{-1}m\} (P_{-2}\otimes P_{-2}).
\end{align*}
The collection of all such vectors $x\otimes x$ comprise the typical fiber.
The $(x\otimes x)(P_{+2}\otimes P_{+2})$ portion of $x\otimes x$ is
identified with the tangent vector and its norm. The $(x\otimes
x)(P_{-2}\otimes P_{-2})$ portion is identified with the momentum vector and
its norm. The 
\begin{equation*}
(x\otimes x)(P_{+2}\otimes P_{-2}+P_{-2}\otimes P_{+2})
\end{equation*}
subspace is associated with the flow of action and its norm. The condition
of minimum action translates into the condition that $x$ has the form 
\begin{equation*}
x=\frac{1}{x_{0}}(x_{0}+x_{1}\mathbf{e}_{1})(x_{0}+x_{2}\mathbf{e}_{2})\; 
\text{for}\; x_{0}\neq 0.
\end{equation*}
It remains to determine how the fibers at different locations are connected.
The description of this connection should also lead to a description of
extended motions on this manifold.


\vskip1pc {\obeylines
\noindent }

{Robert W. Johnson \noindent }

{878 Sunnyhills Road }

{\noindent Oakland, CA 94610 \noindent }


\begin{thebibliography}{9}
\bibitem{1}  {\small R.W. Johnson, Found. Phys. V26(2), (1996) 197. }

\bibitem{2}  {\small N. Salingaros, J. Math. Phys. 25(4), (1984) 738. }

\bibitem{3}  {\small D. Hestenes, \textit{Space-Time Algebra} (Gordon \&
Breach, New York, 1966); P. Lounesto, Found. Phys. \textbf{23}(9), (1993)
1203; K. R. Greider, Found. Phys. \textbf{14}(6), (1984) 467. }

\bibitem{4}  {\small V. I. Arnold, Mathematical Methods of Classical
Mechanics, 2nd Ed. (Springer-Verlag, New York, 1989) p. 253. }

\bibitem{5}  {\small P. J. E. Peebles, Principles of Physical Cosmology,
(Princeton University Press, Princeton, 1993) p. 244. }
\end{thebibliography}
\end{document}